\documentclass[12pt]{article}
\usepackage{fullpage,graphicx,graphics,psfrag,amsmath,amsfonts,verbatim,siunitx,url,csquotes}
\usepackage[small,bf]{caption}
\usepackage{amssymb,enumitem}
\usepackage{booktabs}

\newcommand{\BEAS}{\begin{eqnarray}}
\newcommand{\EEAS}{\end{eqnarray}}
\newcommand{\BEA}{\begin{eqnarray}}
\newcommand{\EEA}{\end{eqnarray}}
\newcommand{\BEQ}{\begin{equation}}
\newcommand{\EEQ}{\end{equation}}
\newcommand{\BIT}{\begin{itemize}}
\newcommand{\EIT}{\end{itemize}}
\newcommand{\BNUM}{\begin{enumerate}}
\newcommand{\ENUM}{\end{enumerate}}

\newcommand{\BA}{\begin{array}}
\newcommand{\EA}{\end{array}}

\newcommand{\eg}{{\it e.g.}}
\newcommand{\ie}{{\it i.e.}}

\newcommand{\reals}{{\mbox{\bf R}}}

\newcommand{\symm}{{\mbox{\bf S}}}  

\newcommand{\diag}{\mathop{\bf diag}}

\newcommand{\argmin}{\mathop{\rm argmin}}

\newcounter{algorithmctr}[section]
\renewcommand{\thealgorithmctr}{\thesection.\arabic{algorithmctr}}
\newenvironment{algdesc}%
   {\refstepcounter{algorithmctr}\begin{list}{}{%
       \setlength{\rightmargin}{0.02\linewidth}%
       \setlength{\leftmargin}{0.02\linewidth}}%
       \rmfamily\small
       \item[]{\setlength{\parskip}{0ex}\hrulefill\par%
        \nopagebreak{\bfseries\textsf{Algorithm \thealgorithmctr~}}}}%
   {{\setlength{\parskip}{-1ex}\nopagebreak\par\hrulefill} \end{list}}

\title{Fitting a Kalman Smoother to Data}
\author{Shane Barratt and Stephen Boyd}

\begin{document}
\maketitle

\begin{abstract}
This paper considers the problem of fitting
the parameters of a Kalman smoother to data.
We formulate the Kalman smoothing problem
with missing measurements
as a constrained least squares
problem and provide an efficient method to solve it
based on sparse linear algebra.
We then introduce the Kalman smoother tuning
problem, which seeks to find parameters
that achieve low prediction error on held out
measurements.
We derive a Kalman smoother auto-tuning
algorithm, which is based on the proximal gradient
method, that finds good, if not the best,
parameters for a given dataset.
Central to our method is the computation of the gradient
of the prediction error with respect to the parameters
of the Kalman smoother; we describe how to compute this at little to no
additional cost.
We demonstrate the method on population migration within the United States
as well as data collected from an IMU+GPS system while driving.
The paper is accompanied by an open-source implementation.
\end{abstract}

\section{Introduction}
Kalman smoothers are widely used to estimate the state of a linear
dynamical system from noisy measurements \cite{swerling1958proposed,kalman1960new}.
In the traditional formulation, the dynamics and output matrices are
considered fixed attributes of the system; the covariance matrices
of the process and sensor noise are tuned by the designer,
within some limits, to obtain good performance in simulation or on
the actual system.  For example, it is common to use noise levels
in the Kalman smoother well in excess of the actual noise, to
obtain practical robustness \cite[\S 8]{gelb1974applied}.

In this paper we take a machine learning approach to the problem of
tuning a Kalman smoother.  We start with the observation that (by
our definition) only the output is observed.  This implies
that the only way we can verify that a Kalman smoother is working
well is to compare the outputs we predict with those that
actually occur, on new or unseen test data, \ie, data that was not used
by the Kalman smoother.
In machine learning terms, we would consider this output prediction
error to be our error, with the goal of minimizing it.
We consider the noise covariance matrices,
as well as the system matrices, as \emph{parameters} that
can be varied
to obtain different estimators, in this case, different Kalman
smoothers.
These are varied, within limits, to obtain good test performance; this
final Kalman smoother can then be checked on entirely new data.

To do this we formulate the Kalman smoothing problem, with missing 
observations, as a simple least squares problem, with a coefficient
matrix that depends on the parameters, \ie, the system
and noise covariance matrices.  We show how to efficiently compute
the derivative of the test error with respect to the parameters,
and use a simple proximal gradient method to update them to improve
the test error.
This method yields a Kalman smoother \emph{auto-tuning} method.
It uses one or more observed output sequences, and the usual
prior knowledge in determining the starting system matrices as well
as a description of the set over which we are allowed to vary them.

The same formulation works for tuning robust Kalman smoothers, where
the process and sensor noises are assumed to have a non-Gaussian
distribution, typically with fatter tails.  In this case the
least squares formulation of the Kalman smoother becomes a convex
optimization problem, and the effect of the parameters
is even less obvious, and therefore harder to tune manually.
Our auto-tuning method extends immediately to such problems.

In summary, the contributions of this paper are:
\begin{itemize}
  \item We describe a Kalman smoother auto-tuning method that requires
  only a dataset of measurements, which may have missing entries.
  \item We describe an efficient method
  for computing the gradient of the prediction error with
  respect to the Kalman smoother parameters, that incurs little to no additional
  computational cost on top of already smoothing.
  \item We provide an open-source implementation of the aforementioned ideas and illustrate the
  method via numerical examples that use real data.
\end{itemize}

\section{Related work}
The Kalman filter was independently invented by
Swerling \cite{swerling1958proposed} and Kalman \cite{kalman1960new}
around 1960, and one of its original applications
was for space aircraft tracking
in the Apollo navigation system \cite{mcgee1985discovery}.
The Kalman filter assumes \emph{a priori} knowledge of the system matrices
and noise statistics. 
Indeed, in his ground-breaking paper, Kalman remarked on the difficulty of identifying
such parameters:
\begin{displayquote}
In real life, however, the situation is usually reversed. One is
given the covariance matrix [of the state] and the problem
is to get [the dynamics] and the statistical properties of [the disturbance]. This is a subtle
and presently largely unsolved problem in experimentation and
data reduction.
\end{displayquote}

Despite its wide use and success, practitioners employing the Kalman smoother
still have to resort to manually tuning its parameters.
As a result, many have proposed methods for instead
automatically tuning the parameters in Kalman smoothers.
One of the first methods proposed was to jointly learn
the parameters and state/output sequence using 
expectation-maximization \cite{shumway1982approach}.
More recent approaches employ different optimization approaches, including
the simplex algorithm \cite{powell2002automated}, coordinate descent \cite{abbeel2005discriminative},
genetic algorithms \cite{oshman2000optimal}, nonlinear programming using finite differencing
to estimate the gradient \cite{asmar2012nonlinear}, Bayesian optimization \cite{chen2018weak},
and reinforcement learning \cite{goodall2007intelligent}.

Our approach is inspired
by previous research on automatically tuning hyper-parameters in least
squares \cite{barratt2019least}.
Our paper departs from prior work on tuning Kalman filters in several ways.
Since our Kalman smoother can deal with missing measurements, we can
hold out measurements and use those to evaluate the smoother.
Also, our method makes explicit use of the gradient of the loss with
respect to the parameters, leading to a more efficient optimization
algorithm.

\section{Kalman smoother}

\paragraph{System model.}
We consider a linear system with dynamics
\BEQ\label{e-dynamics}
x_{t+1} = Ax_t + w_t, \quad t=1,\ldots,T-1,
\EEQ
and output or sensor measurements
\BEQ\label{e-sensor}
y_t = Cx_t + v_t, \quad t=1,\ldots,T.
\EEQ
Here
$x_t\in\reals^n$ is the state,
$w_t\in\reals^n$ is the process noise,
$y_t\in\reals^p$ is the output or sensor measurement,
and $v_t\in\reals^p$ is the sensor noise, at time $t$.
The matrix $A\in\reals^{n\times n}$ is the state dynamics matrix
and $C\in\reals^{p \times n}$ is the output matrix.

We make the standard statistical assumptions that $w_1, \ldots, w_{T-1}$
are IID $\mathcal N(0,W)$
and $v_1, \ldots, v_{T}$ are IID $\mathcal N(0,V)$,
where the symmetric positive definite matrices $W$ and $V$ are the process and sensor 
noise covariance matrices, respectively.

\paragraph{Missing measurements.}
We assume throughout that only the sequence $y_t$ is observed.  Indeed,
we will assume that not all of the measurements are available to us.
To model this, we modify the output equation \eqref{e-sensor} so
that $y_t \in (\reals \cup \{?\})^p$, where ? denotes a missing value.
We have
\BEQ\label{e-sensor-missing}
(y_t)_i = (Cx_t + v_t)_i, \quad (t,i) \in \mathcal K,
\EEQ 
where $\mathcal K \subseteq \{1,\ldots, T\} \times \{1,\ldots, p\}$
is the set of (scalar) outputs that are available.
For $(t,i) \not\in \mathcal K$, we take $(y_t)_i = \; ?$.
We refer to entries of $y_t$ that are real as \emph{known measurements}
and the entries of $y_t$ that have the value ? as \emph{missing 
measurements}.

\paragraph{Smoothing.}
The goal in smoothing is to reconstruct or approximate 
the missing measurements given the known measurements.
Since the outputs and states are jointly Gaussian, the maximum
likelihood and conditional mean estimates of the missing output
values are the same, and can be found as the solution of the constrained
least squares problem
\BEQ
\begin{array}{ll}
\mbox{minimize} &
\sum_{t=1}^{T-1} \|W^{-1/2} (\hat{x}_{t+1} - A\hat{x_t})\|_2^2 + \\[.2cm]
& \sum_{t=1}^T \|V^{-1/2} (\hat{y_t} - C\hat{x_t})\|_2^2, \\[.2cm]
\mbox{subject to} & (\hat y_t)_i = (y_t)_i, \quad (i,t) \in \mathcal K,
\end{array}
\label{eq:kf}
\EEQ
with variables $\hat x_1,\ldots,\hat x_T$ and $\hat y_1,\ldots,\hat y_T$.

Also, the problem has a simple and widely used recursive solution for $\hat x_t$
when $\mathcal K = (1,\ldots,m) \times (1,\ldots,T)$, and also when $T\rightarrow\infty$.
This recursive solution is often referred to as the Kalman filter \cite{kalman1960new}.

\paragraph{Solving the Kalman smoothing problem.}
There are many ways to solve the Kalman smoothing problem \eqref{eq:kf}.
One method is to eliminate the equality constraint [\S 4.2.4]\cite{boyd2004convex}
and solve the resulting unconstrained least squares problem,
which has a banded coefficient matrix.
This method has time and space complexity of order $T(n+p)^2$.
We give some details on another method that has roughly the same complexity,
but is simpler since it does not require eliminating the equality constraints.

Let $N=T(n+p)$ and define the vector $z\in\reals^{N}$ as
$z = (\hat{x}_1,\ldots,\hat{x}_T,\hat{y}_1,\ldots,\hat{y}_T)$.
Using the variable $z$, we can express the estimation problem \eqref{eq:kf} compactly as
the constrained least squares problem
\BEQ
\begin{array}{ll}
\mbox{minimize} &
\frac{1}{2} \|Dz\|_2^2 \\[.05cm]
\mbox{subject to} & Bz=c,
\end{array}
\label{eq:lstsqkf}
\EEQ
where $B\in\reals^{|\mathcal K| \times N}$ is a selector matrix and $c\in\reals^{|\mathcal K|}$ contains
the corresponding entries of $y_t$. 
Concretely, if we assume that $\mathcal K$ is ordered, then if $\mathcal K_j=(i, t)$, the $j$th row
of $B$ is $e_{Tn+tp+i}$ and the $j$th entry of $c$ is $(y_t)_i$.
The matrix $D\in\reals^{N-n \times N}$ is given by
\[
D = \left[
\begin{array}{cc}
D_{11} & 0 \\
D_{21} & D_{22}
\end{array}
\right],
\]
where
\[
\begin{array}{lll}
D_{11} &=&
\hspace*{-.15cm}
\left[
\begin{array}{cccc}
-W^{-1/2}A & W^{-1/2} & 0 & 0\\
0 &  \ddots & \ddots & 0\\
0 & 0 & -W^{-1/2}A & W^{-1/2} \\
\end{array}
\right], \\[.7cm]
D_{21} &=&
\hspace*{-.15cm}
\left[
\begin{array}{ccc}
-V^{-1/2}C & 0 & 0 \\
0 & \ddots & 0 \\
0 & 0 & -V^{-1/2}C
\end{array}
\right], \\[.7cm]
D_{22} &=&
\hspace*{-.15cm}
\left[
\begin{array}{ccc}
V^{-1/2} & 0 & 0 \\
0 & \ddots & 0 \\
0 & 0 & V^{-1/2} \\
\end{array}
\right]. \\
\end{array}
\]
The matrices $D$ and $B$ are evidently very sparse,
since each have a density of approximately $\frac{1}{N}$.

The optimality conditions for \eqref{eq:lstsqkf} can be expressed as
\[
\left[
\begin{array}{ccc}
0 & D^T & B^T \\
D & -I & 0 \\
B & 0 & 0
\end{array}
\right]
\left[
\begin{array}{c}
z \\ v \\ \eta
\end{array}
\right]
=
\left[
\begin{array}{c}
0 \\ 0 \\ c
\end{array}
\right],
\]
where $\eta\in\reals^{|\mathcal K|}$ is the dual variable for the equality constraint
and $v=Dz$.
The KKT matrix, denoted by
\[
M=\left[
\begin{array}{ccc}
0 & D^T & B^T \\
D & -I & 0 \\
B & 0 & 0
\end{array}
\right],
\]
is also very sparse, since $B$ and $D$ are sparse.

We assume for the remainder of the paper that
$M$ is full rank (if it is not, we can add a small amount of regularization to make it invertible).
Therefore we can solve the KKT system using any method for solving a sparse
system of linear equations, \eg, a sparse LU factorization \cite{davis2004algorithm}. 
Since the sparsity pattern is banded (when re-ordered the right way),
the complexity of the sparse LU factorization will be linear in $T$.
We have also observed this to be true in practice (see figure~\ref{fig:timings}).

\paragraph{Judging a Kalman smoother.}

Suppose we have gathered a sequence of outputs denoted $y_1,\ldots,y_T \in (\reals \cup \{?\})^p$.
We can judge how well a Kalman smoother is working on this sequence of observations by
obscuring a fraction of the known outputs and comparing
the outputs predicted by the Kalman smoother to those that actually occurred.

The first step in judging a Kalman smoother is to
mask some fraction (\eg, 20\%) of the non-missing entries in the observations, denoted by the set
$\mathcal M_i \subseteq (1,\ldots,T) \times (1,\ldots,m)$,
resulting in a \emph{masked trajectory} $\tilde{y}_1, \ldots, \tilde{y}_T$.
That is, we let $(\tilde{y}_t)_i=\;?$ for $(i,t)\in\mathcal M$ and $(\tilde y_t)_i = (y_t)_i$ for $(i, t) \not \in \mathcal M$.

We then solve the smoothing problem \eqref{eq:kf} with $y_t=\tilde y_t$ and known set $\mathcal K \setminus \mathcal M$,
resulting in a predicted output trajectory $\hat y_1,\ldots, \hat y_T$.

In order to judge the Kalman smoother, we calculate the squared difference between the predicted output trajectory
and the actual trajectory
in the entries that we masked, which is given by
\BEQ
L = \sum_{(i,t) \in \mathcal M} \left((\hat y_t)_i - (y_t)_i\right)^2.
\label{eq:error}
\EEQ
We refer to this quantity as the \emph{prediction error}; the goal in the sequel will be to adjust
the parameters to minimize this error.
We note that the entries in the output should
been suitably scaled or normalized such that \eqref{eq:error}
is a good measure of prediction error for the given application.

\section{Kalman smoother auto-tuning}
In this section we describe how to automatically tune the parameters in a Kalman smoother (that is,
the dynamic matrices and covariance matrices) to minimize the prediction error on
the held-out measurements \eqref{eq:error}.
Once the parameters have been tuned, the Kalman smoother can be tested on another (unseen)
output sequence.

\begin{figure*}[t]
\begin{algdesc}
\label{alg:ksat}
\emph{Kalman smoother auto-tuning.}
\begin{tabbing}
    {\bf given} initial hyper-parameter vector $\theta^1 \in \Theta$,
    initial step size $t^1$, number of iterations $n_\mathrm{iter}$,\\
    \qquad \=\ tolerance $\epsilon$.\\
    {\bf for} $k=1,\ldots,n_\mathrm{iter}$\\
      \qquad \=\ 1.\ \emph{Filter the output sequence}. Let
      $\hat y_1,\ldots,\hat y_T$ be the solution to \eqref{eq:kf}.\\
      \qquad \=\ 2.\ \emph{Compute the gradient of the prediction error}.
    $g^k = \nabla_\theta L(\theta)$.\\
      \qquad \=\ 3.\ \emph{Compute the gradient step}. $\theta^{k+1/2} = \theta^k - t^k g^k$.\\
      \qquad \=\ 4.\ \emph{Compute the proximal operator}. 
    $\theta^\mathrm{tent} = \mathrm{\bf prox}_{t^k r} (\theta^{k+1/2})$.\\
      \qquad \=\ 5. {\bf if} $F(\theta^\mathrm{tent}) \leq F(\theta^k)$:\\
      \qquad \qquad \=\ \emph{Increase step size and accept update.}
      $t^{k+1} = (1.5)t^k; \quad \theta^{k+1} = \theta^\mathrm{tent}$.\\
      \qquad \qquad \=\ \emph{Stopping criterion.} {\bf quit} if
      $\|(\theta^k - \theta^{k+1})/t^k + (g^{k+1} - g^k) \|_2 \leq \epsilon$.\\
      \qquad \=\ 6. {\bf else} \emph{Decrease step size and reject update}.
      $t^{k+1} = (0.5)t^k; \quad \theta^{k+1} = \theta^k$.\\
    {\bf end for}
\end{tabbing}
\end{algdesc}
\vspace{-1em}
\end{figure*}

\paragraph{Kalman smoother parameters.}
A Kalman smoother has four parameters, which we denote by
\[
\theta=(A,W^{-1/2},C,V^{-1/2})\in
\reals^{n\times n}\times\reals^{n\times n}\times\reals^{p\times n}\times\reals^{p \times p}.
\]
Evidently, this parametrization of the Kalman smoother is not unique.
For example, if $T\in\reals^{n\times n}$ is invertible, then
$\tilde x_t = T x_t$, $\tilde A = TAT^{-1}$,
$\tilde W=T^{-1}WT^{-T}$, $\tilde C=CT^{-1}$,
and $\tilde V=V$ gives another representation of \eqref{eq:kf}.
As another example, scaling $W$ and $V$ by $\alpha > 0$
gives an equivalent representation of \eqref{eq:kf}.

\subsection{Auto-tuning problem}
The prediction error $L$ in \eqref{eq:error} is a function of the parameters,
and from here onwards we denote that function by $L(\theta)$.
To tune the Kalman smoother, we propose solving the optimization problem
\BEQ
\begin{array}{ll}
\mbox{minimize} &
F(\theta) = L(\theta) + r(\theta), \\
\end{array}
\label{eq:tune}
\EEQ
with variable $\theta$ (the parameters of the Kalman smoother),
where $r:\Theta \to \reals$ is a regularization function.
Here $\Theta$ denotes the set of allowable parameters and can, for example,
include constraints on what parameters we are allowed to change.
(The function $r$ evaluates to $+\infty$ for $\theta \not \in \Theta$,
thus constraining $\theta$ to be in $\Theta$.)

The objective function $F:\Theta \to \reals$ is composed of two parts: the prediction
error and the regularization function.
The first term here encourages the Kalman filter to have the same outputs
as those observed,
and the second term encourages the parameters to be simpler
or closer to an initial guess.

\paragraph{Regularization functions.}
There are many possibilities for the regularization function $r$;
here we describe a few.
Suppose we have some initial guess for $A$, denoted
$A_\mathrm{nom}$.
We could then penalize deviations of $A$ from $A_\mathrm{nom}$ by letting, \eg,
\[
r(\theta) = \|A-A_\mathrm{nom}\|_F^2.
\]

As another example, suppose we suspected that $C$ was low rank; then we could use
\[
r(\theta) = \|C\|_*,
\]
where $\|C\|_*$ is the nuclear norm of $C$, \ie,
the sum of the singular values of $C$.
This regularizer encourages $C$ to be low rank.
Of course, any combination of these regularization functions is possible.

\paragraph{Allowable sets.}
There are also many possibilities for $\Theta$, the allowable set of parameters.
One option is to only allow certain entries of $A$ to vary by letting the set of allowable $A$ matrices
be
\[
\{A \mid A_{ij} = (A_\mathrm{nom})_{ij}, (i,j) \in \Omega\}
\]
for some set $\Omega$.
If we wanted to keep $A$ fixed, we could let $\Theta=\{A_\mathrm{nom}\}$.
Another sensible option is to let $A$ vary within a box by letting
the set of allowable $A$ matrices be
\[
\{A \mid \|A - A_\mathrm{nom}\|_\infty \leq \rho\},
\]
for some nominal guess $A_\mathrm{nom}$ and hyper-parameter $\rho > 0$.

\subsection{Solution method}
\label{sec:solution_method}

The auto-tuning problem \eqref{eq:tune} is in general nonconvex, even if $\Theta$ and $r$ are convex,
so it is very difficult to solve exactly.
Therefore, we must resort to a local or heuristic optimization method to (approximately) solve it. 
There are many methods that we could use
to (approximately) solve the auto-tuning problem (see, \eg,
\cite{douglas1956numerical,lions1979splitting,shor1985minimization,nesterov2013introductory}).
In this paper we employ one of the simplest, the proximal gradient method \cite{martinet1970breve,nesterov2013gradient},
since $F$ is differentiable in $\theta$ (see below).

The proximal gradient method is described by the iteration
\[
\theta^{k+1} =
\mathrm{\bf prox}_{t^k r}
(\theta^k - t^k \nabla_\theta L(\theta)),
\]
where $k$ is the iteration number, $t^k > 0$ is a step size, and the proximal operator
of $tr(\cdot)$ is defined as
\[
\mathrm{\bf prox}_{tr}(\nu) = \argmin_{\theta \in \Theta}
\left(t r(\theta) + (1/2)\|\theta-\nu\|_2^2\right).
\]
When $\Theta$ is a convex set and $r$ is convex, evaluating the proximal operator
of $r$ requires solving a (small) convex optimization problem.
Also, the proximal operator often has a (simple) closed-form expression \cite{parikh2014proximal}.
We note that $r$ need not be differentiable.

We employ the proximal gradient method
with the adaptive step size scheme and stopping condition described in \cite{barratt2019least}.
The full algorithm for Kalman smoother auto-tuning is summarized in Algorithm \ref{alg:ksat}.

\paragraph{Computing the gradient.}

Evidently, the proximal gradient method requires computing the gradient of the prediction error
with respect to the parameters,
denoted $\nabla_\theta L(\theta)$.
The sensitivity analysis of Kalman smoothing has previously been considered in
the forward direction, \ie, how changes in the parameters affect the output \cite[\S 7]{gelb1974applied}.
Justification for our derivation can be found in \cite[\S3.5]{barratt2019least}.

To do this, we first form the gradient of $L$ with respect to $\hat{y}_1,\ldots,\hat{y}_T$,
given by
\[
\nabla_{(\hat{y}_t)_i} L = \begin{cases}
2\left((\hat{y}_t)_i - (y_t)_i\right) & (i,t) \in \mathcal M, \\
0 & \text{otherwise}.
\end{cases}
\]
Next we form the gradient of $L$ with respect to the solution to \eqref{eq:lstsqkf},
which is given by 
\[
g = \nabla_{(z,\eta,v)} L=(0,\nabla_{y_1} L, \ldots, \nabla_{y_T} L, 0, 0).
\]
Since
\[
(z, \nu, \eta) = M^{-1} 
\]
Next we solve the linear system
\[
M
\left[\begin{array}{c}
q_1 \\ q_2 \\ q_3 \end{array}\right] = -g,
\]
which only requires a backsolve if we have already factorized $M$.
Since the KKT system $M$ is invertible, the prediction error is indeed differentiable.

The next step is to form the gradient of $L$ with respect to the coefficient matrix $D$, which
is given by
\[
G = \nabla_D L = D(q_1z^T + zq_1^T).
\]
Since $\theta$ only affects $D$ at certain entries,
we only need to compute $G$ at those entries.
That is, we compute $G$ at the entries
\[
G = \left[
\begin{array}{cc}
G_{11} & 0 \\
G_{21} & G_{22}
\end{array}
\right],
\]
where
\[
\begin{array}{lll}
G_{11} &=&
\left[
\begin{array}{cccc}
G_{11} & G_{12} & 0 & 0\\
0 &  \ddots & \ddots & 0\\
0 & 0 & G_{(T-1)(T-1)} & G_{(T-1)T} \\
\end{array}
\right], \\[.7cm]
G_{21} &=&
\left[
\begin{array}{ccc}
G_{T1} & 0 & 0 \\
0 & \ddots & 0 \\
0 & 0 & G_{(2T-1)T}
\end{array}
\right], \\[.7cm]
G_{22} &=&
\left[
\begin{array}{ccc}
G_{(T+1)T} & 0 & 0 \\
0 & \ddots & 0 \\
0 & 0 & G_{(2T-1)(2T-1)} \\
\end{array}
\right], \\
\end{array}
\]
which we can efficiently do since
\[
G_{ij} = (Dq_1)_i z_j + (Dz)_i (q_1)_j.
\]

The final step is to form the gradients with respect to the parameters,
which are given by
\[
\begin{array}{ccl}
\nabla_A L &=& -(W^{-1/2})^T \sum_{t=1}^{T-1}G_{tt},\\[.1cm]
\nabla_{W^{-1/2}} L &=& \sum_{t=1}^{T-1} G_{t(t+1)} - G_{tt}A^T,\\[.1cm]
\nabla_C L &=& -(V^{-1/2})^T \sum_{t=1}^T G_{(T-1+t)t},\\[.1cm]
\nabla_{V^{-1/2}} L &=& \sum_{t=1}^T G_{(T+t)(T+1+t)} - G_{(T-1+t)t}C^T.
\end{array}
\]
The complexity of computing the gradient is roughly the same
complexity as solving the original problem,
since it requires the solution of another linear system.
However, the time required to
compute the gradient is often lower since we cache the factorization
of $M$.

\section{Experiments}

In this section, we describe our implementation of Kalman smoother
auto-tuning, as well as the results
of some numerical experiments that illustrate the method.
All experiments were performed on a single core of
an unloaded Intel i7-8770K CPU.

\paragraph{Reference implementation.}

We have implemented the Kalman smoother auto-tuning method
described in this paper as an open-source Python package,
available at
\[
\verb|https://github.com/cvxgrp/auto_ks|.
\]
Our CPU-based implementation
has methods for performing Kalman smoothing with missing measurements
and for tuning the matrices in the Kalman smoother (Algorithm \ref{alg:ksat}).
Our only dependencies are \verb|scipy| \cite{scipy}, which we use for sparse linear algebra,
and \verb|numpy| \cite{walt2011numpy}, which we use for dense linear algebra.

\begin{figure}
    \centering
    \includegraphics[width=0.5\linewidth]{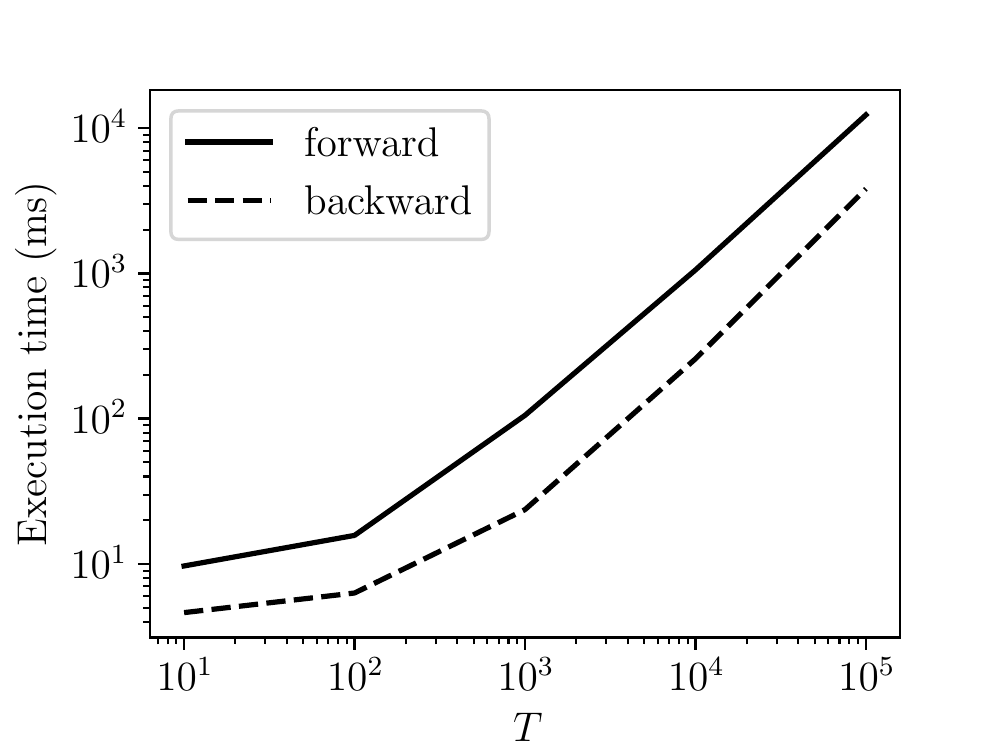}
    \caption{Method timings for a random problem with $n=p=10$.}
    \label{fig:timings}
    \vspace{-1em}
\end{figure}

\paragraph{Performance.}
We ran our Kalman smoothing function on random problems with $n=p=10$.
Fig. \ref{fig:timings} shows the execution time,
averaged over ten runs, of solving the smoothing problem (denoted as \emph{forward} in the figure),
as well as computing the derivative with respect to the the parameters (denoted as \emph{backward} in the figure).
As expected, the time required to compute the solution and its derivative
is roughly linear in the length of the sequence $T$.
Empirically, we found that the time required to compute the derivative is roughly
half of the time required to compute the solution.
We remark that our method is very efficient and effortlessly
scales to extremely large problem sizes.

\subsection{Human migration example}
Suppose we have $n$ states, where the $i$th state has a population $x_i$.
At some cadence, say yearly, a fraction of people in each state decide to
move to another state.
We take noisy measurements of the population in some of the states
and wish to infer the population in every state, including those we have not
even measured.

\paragraph{Dynamics.}
The vector $x_t\in\reals^n$ denotes the population in each state at year $t$.
The dynamics are described by
\[
x_{t+1} = A x_t + w_t,
\]
where $A\in\reals_+^{n\times n}$.
Here $A_{ij}$ denotes the fraction of the population in state $j$
that move to state $i$ each year.

\paragraph{Outputs.}
Each year, we take noisy measurements of the populations in some of the states.
The outputs are described by
\[
y_t = x_t + v_t.
\]
We use the set of known measurements $\mathcal K$ to denote the measurements
we actually have access to.

\paragraph{United States population data.}
We gathered yearly population data (in millions of people)
for the $n=48$ states in the continental U.S
from the U.S. Census Bureau \cite{censusdata}.
The data includes all years from 1900 to 2018.

\paragraph{Experiment details.}
Our goal is to learn the dynamics matrix,
dynamics covariance, and output covariance via Kalman smoother auto-tuning.
To this end, we use the regularization function $r(\theta)=0$ and allowable set
\[
\Theta = \{(A, \diag(w), C_0, \diag(v)) \mid A \in \reals_+^{n \times n}, w,v \in \reals_+^n\}.
\]
We initialize the parameters as
\[
A_0=I, \quad W^{-1/2}_0 = (30)I, \quad C_0=I, \quad V^{-1/2}_0 = (10)I
\]
For each year, we pick 30 out of the 48 states at random to be measured.
In each year, of those measured, we pick 12 at random to be missing
and 5 at random to be part of the test set.
We ran the method for 50 iterations with $t_0=\num{1e-4}$.

\paragraph{Results.}
The prediction error decreased from 0.0097 to 0.0058.
The test error decreased from 0.0041 to 0.0030.
The algorithm took 31 seconds to run.
Besides the purely numerical results, there are interesting interpretations
of the resulting parameters.
For example, we can interpret the off-diagonal entries in the $A$ matrix
as the fraction of the population in one state that migrates to another state
over the course of one calendar year.
The biggest such entries are displayed in Tab. \ref{tab:migration}.

\begin{table}
  \caption{Entries in the learned $A$ matrix.}
  \vspace{1em}
  \centering
  \begin{tabular}{lll}
    \toprule
    From & To & Fraction \\
    \midrule
    Texas & California & 0.56\% \\
    Florida & California & 0.42\% \\
    California & Texas & 0.30\% \\
    Texas & Florida & 0.29\% \\
    Pennsylvania & California & 0.28\% \\
    \bottomrule
  \end{tabular}
  \label{tab:migration}
\end{table}

\subsection{Vehicle smoothing example}
In vehicle smoothing, we have noisy measurements
of the position, velocity, and acceleration of a vehicle over time,
and wish to infer the true position, velocity, and acceleration at
each time step.

\paragraph{Dynamics.}
\label{sec:veh_dynamics}
The state $x_t=(p_t, v_t, a_t)$ is composed of the position $p_t\in\reals^3$,
the velocity $v_t\in\reals^3$, and the acceleration $a_t\in\reals^3$.
The dynamics are described by
\[
x_{t+1} = \begin{bmatrix}
I & hI & 0 \\
0 & I & hI \\
0 & 0 & I
\end{bmatrix} x_t + w_t.
\]
(This system is often referred to as a double integrator,
since the derivative of $p_t$ is $v_t$ and the derivative of $v_t$ is $a_t$.)

\paragraph{Outputs.}
\label{sec:veh_outputs}
Any output vector and linear output matrix is possible.
In this specific example we use
$y_t=(\hat p_t, \hat a_t, (\hat v_t)_1, (\hat v_t)_2)$, so $p=8$.
The output is described by
\[
y_t = \begin{bmatrix}
I & 0 & 0 & 0 \\
0 & 0 & 0 & I \\
0 & I_2 & 0 & 0
\end{bmatrix} x_t + v_t.
\]

\paragraph{Data.}
We used the Sensor Play data recorder iOS application \cite{sensorplay}
to record the acceleration, attitude, latitude, longitude, heading, speed, and altitude
of an iPhone mounted on a passenger vehicle.
We converted the latitude and longitude into local North-East-Up coordinates,
used the heading to convert the speed into a velocity in local coordinates,
and used the attitude to orient the acceleration to local coordinates.
We recorded data for a total of 330 seconds with a sampling frequency of \SI{100}{\hertz},
resulting in $T=33000$ measurements.

\paragraph{Experiment details.}
Our goal is to learn the state and observation covariance matrices,
via Kalman smoother auto-tuning.
To this end, we penalize the off-diagonal entries of 
\[r(\theta)=\alpha\sum_{i\neq j}(W^{-1/2}_{ij})^2 + \alpha\sum_{i\neq j}(V^{-1/2}_{ij})^2,\]
where $\alpha$ is a hyper-parameter (we use $\alpha=\num{1e-4}$) and
\[
\Theta = \{(A_0, W^{-1/2}, C_0, V^{-1/2}) \mid W \in \symm_+^9, V \in \symm_+^8 \}.
\]
We initialize the covariances as
$W^{-1/2}_0 = I$, $V^{-1/2}_0 = (0.01)I$,
and initialize $A_0$ and $C_0$ as given in Sec. \ref{sec:veh_dynamics} and Sec. \ref{sec:veh_outputs} respectively.
We consider all measurement indices where the GPS or velocity change as known
(since GPS is only useful when it changes) and all acceleration indices as known.
We use 20\% of the known position measurements as the missing measurements and another 20\%
as the test measurements. 
We ran the method for 25 iterations with $t_0=\num{1e-2}$.

\begin{figure}
    \centering
    \hspace*{.8em}
    \includegraphics[width=0.5\linewidth]{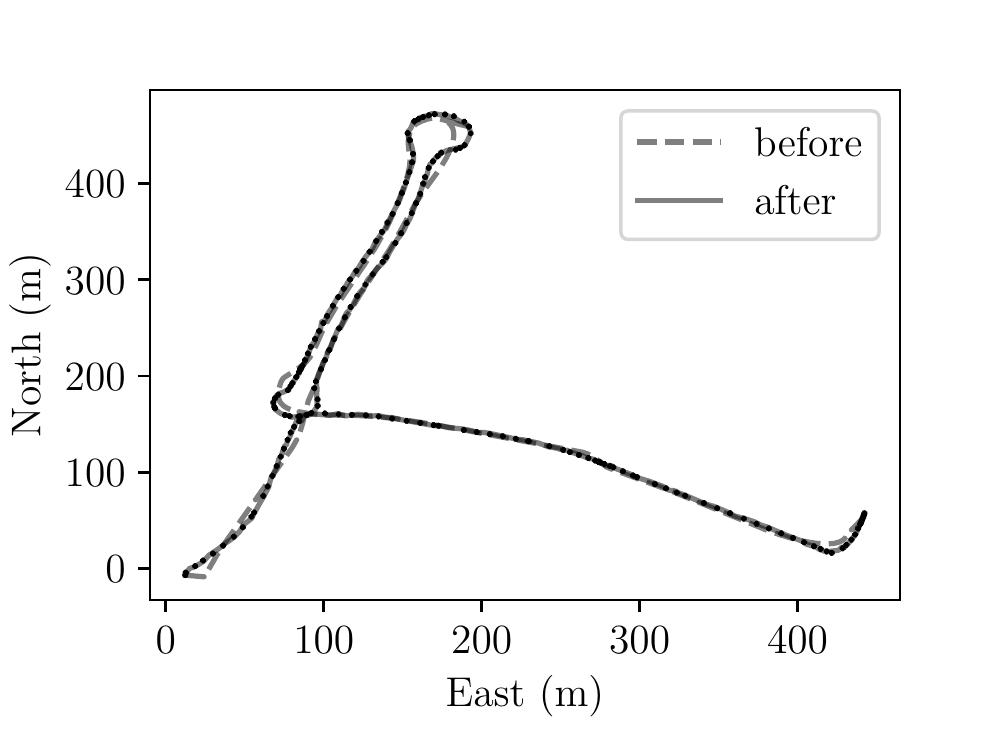}
    \caption{Smoothed position estimates before and after tuning.}
    \label{fig:driving}
\end{figure}  

\paragraph{Results.}
The prediction error decreased from 13.23 to 2.97.
The test error decreased from 16.57 to 1.37.
The algorithm took 135 seconds to run.
The diagonals of the final state and output covariance matrices were 
\[
\begin{array}{lll}
W &=& \diag(5.7, 14.9, 1.1, 0.8, 1.3, 1.0, 1.0, 1.0, 1.3), \\
V &=& \diag(0.2, 0.4, 9.6, 2.7, 2.3, 0.2, 1.4, 2.4).
\end{array}
\]
(Note that these matrices can be scaled and the smoothing result
is the same, so only relative magnitude matters.)
We observe that there is
more state noise in north and east dimensions than up,
which makes sense.
Also, there is less state noise in velocity than in position.
We also observe that there is much higher measurement noise for $z$ direction in GPS,
which is true with GPS.
In Fig. \ref{fig:driving} we show the position estimates before
and after tuning.
Visually, we see significant improvement from tuning.

\section*{Acknowledgments}
S. Barratt is supported by the National Science Foundation Graduate Research Fellowship
under Grant No. DGE-1656518.

\bibliography{refs}

\end{document}